\documentclass[12pt,reqno]{amsart}

\usepackage{amsmath,amsthm,amsfonts,amscd,amssymb,euscript}
\numberwithin{equation}{section} 
\newcommand{\dbar}{\ensuremath{\bar \partial}}

\newcommand{\ad}{\ensuremath{\bar \partial^{*}  }}
\newcommand{\C}{\ensuremath{{\mathbb C}}}
\newcommand{\R}{\ensuremath{{\mathbb R}}}
\newcommand{\N}{\ensuremath{{\mathbb N}}}

\newcommand{\smooth}{\ensuremath{C^{\infty}}}
\newcommand{\smoothc}{\ensuremath{C_0 ^{\infty}}}
\newcommand{\ra}{\ensuremath{C^{\omega}}}
\newcommand{\I}{\ensuremath{\mathcal I}}

\newcommand{\V}{\ensuremath{\mathcal V}}

\newcommand{\B}{\ensuremath{\mathfrak B}}
\newcommand{\T}{\ensuremath{\mathcal T}}
\newcommand{\ls}{\ensuremath{\mathcal L}}
\newcommand{\mt}{\ensuremath{\mathfrak M}}

\newcommand{\ct}{\ensuremath{\mathfrak C}}
\newcommand{\oka}{\ensuremath{\mathcal O}}
\newcommand{\nx}{\ensuremath{{\mathcal N}_x}}
\newcommand{\ny}{\ensuremath{{\mathcal N}_y}}

\newcommand{\cv}{\ensuremath{\mathcal C}}

\begin{document}
\title{Effective vanishing order of the Levi determinant}
\author{Andreea C. Nicoara}

\address{Department of Mathematics, University of Pennsylvania, 209 South $33^{rd}$ St.,  Philadelphia, PA 19104}

\email{anicoara@math.upenn.edu}

\begin{abstract}
On a smooth domain in $\C^n$ of finite D'Angelo $q$-type at a point, an effective upper bound for the vanishing order of the Levi determinant $\text{coeff}\{\partial r \wedge \dbar r \wedge (\partial \dbar r)^{n-q}\}$ at that point is given in terms of the D'Angelo $q$-type, the dimension of the space $n,$ and $q$ itself. The argument uses 
Catlin's notion of a boundary system as well as techniques pioneered by John D'Angelo.
\end{abstract}

\maketitle

\tableofcontents

\section{Introduction}

In his seminal Acta Mathematica paper of 1979 \cite{kohnacta},
Joseph J. Kohn gave a sufficient condition for the subellipticity
of the $\dbar$-Neumann problem on a pseudoconvex domain in $\C^n$
by introducing subelliptic multipliers for the $\dbar$-Neumann
problem as well as an algorithm on these multipliers whose
termination implies subellipticity. Throughout this paper we shall
refer to the latter as the Kohn algorithm. In \cite{kohnacta} Kohn
defined his subelliptic multipliers as germs of $\smooth$
functions, and the sufficient condition he obtained for
subellipticity applies to $\smooth$ pseudoconvex domains. Because
of the peculiar algebraic properties of the ring of $\smooth$
functions, however, he restricted his study of the termination of
his algorithm to the much better behaved ring of real-analytic
functions $\ra.$ It is for such a pseudoconvex domain $\Omega$ in
$\C^n$ with a real-analytic boundary that he proved the
equivalence of the following three properties:

\begin{enumerate}
\item[(i)] subellipticity of the $\dbar$-Neumann problem for $(p,q)$ forms;
\item[(ii)] termination of the Kohn algorithm on $(p,q)$ forms;
\item[(iii)] finite order of contact of holomorphic varieties of complex dimension $q$ with
the boundary of the domain $\Omega.$
\end{enumerate}

The Kohn algorithm generates an increasing chain of ideals of multipliers, one ideal per step. By its termination it is meant that after a finite number of steps, the ideal of multipliers captures a unit; therefore, the algorithm terminates when the entire ring is produced. The first ideal in the Kohn algorithm on $(p,q)$ forms is the real radical of the ideal generated by the defining function of the domain $r$ and the Levi determinant $\text{coeff}\{\partial r \wedge \dbar r \wedge (\partial \dbar r)^{n-q}\}.$ Since $r$ is by definition identically zero on the boundary of the domain, the function that determines the behavior of the Kohn algorithm is  the Levi determinant $\text{coeff}\{\partial r \wedge \dbar r \wedge (\partial \dbar r)^{n-q}\}.$ It follows that  having an effective upper bound for the vanishing order of the Levi determinant is crucial for understanding the Kohn algorithm on smooth domains, as smooth functions, unlike real analytic ones, can vanish to infinite order at a point and yet not be identically zero in a neighborhood. For $q=1$ the Levi determinant is the determinant of the Levi form, so it  is the author's hope this effective upper bound for its vanishing order at a point may serve in other problems that arise in several complex variables.

Two crucial ingredients come into the computation of this effective bound: techniques of John D'Angelo from his work on the openness of condition (iii), also known in the literature as finite D'Angelo type, and Catlin's notions of boundary system and commutator type from \cite{catlinbdry}. The latter were developed by David Catlin in order to show the equivalence of conditions (i) and (iii) for smooth pseudoconvex domains, which he carried out in a series of very deep papers \cite{catlinnec}, \cite{catlinbdry}, and \cite{catlinsubell}. Under the assumption of pseudoconvexity, Catlin showed that the two notions he defined in \cite{catlinbdry}, the multitype and the commutator multitype, equal each other. Finding an effective bound for the vanishing of the Levi determinant, however, will require truncating the Taylor expansion of the defining function $r,$ a step that loses pseudoconvexity. The argument here will only employ the commutator type. For an investigation of the other notion, the Catlin multitype, in the absence of pseudoconvexity, the interested reader should consult Martin Kol\'a\v{r}'s recent work on understanding classes of hypersurfaces of finite Catlin multitype using Chern-Moser invariant theory in \cite{kolar1} and \cite{kolar2}.

We now state the main result of this paper:

\medskip
\newtheorem{mainthm}{Main Theorem}[section]
\begin{mainthm}
 Let $\Omega$ in $\C^n$ be a domain with $\smooth$
boundary. Let $x_0 \in b \Omega$ be a point on the boundary of
the domain such that the D'Angelo $q$-type, namely the order of contact of holomorphic varieties of complex dimension $q$ with the boundary of $\Omega$ at $x_0,$ is finite and equal to $t.$ \label{maintheorem} Let $\lceil t \rceil$ be the roundup of $t,$ i.e. the lowest integer greater than or equal to $t,$ then at $x_0$ the Levi determinant  $\text{coeff}\{\partial r \wedge \dbar r \wedge (\partial \dbar r)^{n-q}\}$ vanishes to order at most $(\lceil t \rceil -2 )^{n-q}.$
\end{mainthm}

Using algebraic geometric techniques involving the ring of holomorphic functions, Yum-Tong Siu proved in \cite{siunote} that a lower bound for the subelliptic gain $\epsilon$ in the $\dbar$-Neumann problem exists that is polynomial in the D'Angelo type $t$ and the dimension $n$ in the case when the pseudoconvex domain is defined by a function $r(z)$ of the special type $$r(z)=Re \, {z_n}+\sum_{j=1}^N \, |f_j(z_1, \dots, z_{n-1})|^2,$$ where $N \geq n$ and $f_j$ is holomorphic for all $1 \leq j \leq N.$ One of the technical results in his paper involved precisely the computation of the order of vanishing of the Levi determinant (his result only applies to $(0,1)$ forms, so for $q=1$) in terms of $t$ and $n.$
By contrast, we carry out this computation here for any smooth domain and any level of forms $q$ using the D'Angelo-Catlin machinery as standard techniques in algebraic geometry that work in the ring of holomorphic functions do not apply in the ring of smooth functions, a far more ill-behaved ring. The main theorem is the first step in the quest for an effective lower bound for the subelliptic gain $\epsilon$ in the $\dbar$-Neumann problem for a smooth pseudoconvex domain that uses the Kohn algorithm. It should be noted that Catlin already obtained in \cite{catlinsubell} a lower bound $$\epsilon \geq \tau^{-n^2 \, \tau^{n^2}}$$ that holds for any smooth pseudoconvex domain and is exponential in $\tau,$ a finite type notion similar to the D'Angelo type. However, Catlin used his own techniques that are completely unrelated to the Kohn algorithm, and $\tau$ is only equal to the D'Angelo type $t$ if $q=1;$ otherwise, the relationship between $\tau$ and $t$ is not clear. Both of these notions will be discussed in Section~\ref{findap}.

This paper is organized as follows: Section~\ref{kohnalg} is
devoted to the Kohn algorithm. Section~\ref{findap} defines finite D'Angelo type and outlines the properties necessary for the proof of the Main Theorem~\ref{maintheorem}. Section~\ref{mtbssect} introduces Catlin's boundary systems as well as his multitype and commutator multitype. The Main Theorem~\ref{maintheorem} is proven in Section~\ref{mainthmpf}.

The author wishes to thank Charles L. Fefferman for a number of very useful discussions.

\section{The Kohn algorithm}
\label{kohnalg}

It already follows from Joseph J. Kohn's solution to the
$\dbar$-Neumann problem in \cite{dbneumann1} and \cite{dbneumann2}
for strongly pseudoconvex domains as well as from the weighted
estimates for pseudoconvex domains done H\"{o}rmander in
\cite{hormest} and Kohn in \cite{kohnweakpsc} that the
$\dbar$-Neumann problem is elliptic inside the domain $\Omega,$
so the study of subellipticity only needs to be
conducted on the boundary of the domain $b \Omega.$

We start with Joseph J. Kohn's definition of what it
means for the $\dbar$-Neumann problem on $(p,q)$ forms to be
subelliptic followed by his definition of a subelliptic multiplier
from \cite{kohnacta}. We refer the reader to this same paper for
details and motivation regarding the setup of the
$\bar\partial$-Neumann problem:

\medskip
\newtheorem{subellgain}{Definition}[section]
\begin{subellgain}
Let $\Omega$ be a domain in $\C^n$ and let $x_0 \in
\overline{\Omega}.$ The $\dbar$-Neumann problem on $\Omega$ for
$(p,q)$ forms is said to be subelliptic at $x_0$ if there exist a
neighborhood $U$ of $x_0$ and constants $C, \epsilon > 0$ such
that
\begin{equation}
||\varphi\, ||_{\, \epsilon}^2 \leq C \, ( \, ||\,\dbar \,
\varphi\,||^2_{\, 0} + ||\, \ad  \varphi \,||^2_{\, 0} +
||\,\varphi \,||^2_{\, 0} \,) \label{subellestproblem}
\end{equation}
for all $(p,q)$ forms $\varphi \in \smoothc (U) \cap Dom (\ad),$
where $||\, \cdot \, ||_{\, \epsilon}$ is the Sobolev norm of
order $\epsilon$ and $||\, \cdot \,||_{\, 0}$ is the $L^2$ norm.
\label{subellgaindef}
\end{subellgain}

\medskip
\newtheorem{subellmult}[subellgain]{Definition}
\begin{subellmult}
Let $\Omega$ be a domain in $\C^n$ and let $x_0 \in
\overline{\Omega}.$ A $\smooth$ function $f$ is called a
subelliptic multiplier at $x_0$ for the $\bar\partial$-Neumann
problem on $\Omega$ for $(p,q)$ forms if there exist a neighborhood $U$ of $x_0$ and
constants $C, \epsilon > 0$ such that
\begin{equation}
||\, f \varphi\, ||_{\, \epsilon}^2 \leq C \, ( \, ||\,\dbar \,
\varphi\,||^2_{\, 0} + ||\, \ad  \varphi \,||^2_{\, 0} +
||\,\varphi \,||^2_{\, 0} \,) \label{subellest}
\end{equation}
for all $(p,q)$ forms $\varphi \in \smoothc (U) \cap Dom (\ad).$
We will denote by $I^q (x_0)$ the set of all subelliptic
multipliers at $x_0.$ \label{subellmultdef}
\end{subellmult}

\medskip\noindent Several remarks about these two definitions are necessary:

\smallskip \noindent (1) If there exists a subelliptic multiplier $f \in I^q (x_0)$ such that $f(x_0) \neq 0$, then a
subelliptic estimate holds at $x_0$ for the $\dbar$-Neumann
problem.

\smallskip\noindent (2) As explained at the beginning of this section, if $x_0 \in \Omega$ then automatically estimate
\ref{subellest} holds at $x_0$ with the largest possible
$\epsilon$ allowed by the $\dbar$-Neumann problem, namely
$\epsilon = 1.$ This says the problem is elliptic rather than
subelliptic inside.

\smallskip\noindent (3) The previous remark implies that if $x_0
\in b \Omega$ but $f = 0$ on $U \cap b \Omega,$ then estimate
\ref{subellest} again holds for $\epsilon = 1.$ This is the case
if we set $f = r,$ where $r$ is the defining function of the
domain $\Omega.$

\smallskip\noindent (4) If $x_0 \in b \Omega,$ the highest
possible gain in regularity in the $\dbar$-Neumann problem is given
by $\epsilon=\frac{1}{2}$ under the strongest convexity
assumption, namely strong pseudoconvexity of $\Omega$, as proved
by Kohn in \cite{dbneumann1} and \cite{dbneumann2}.

\smallskip\noindent (5) This non-ellipticity of the
$\dbar$-Neumann problem is coming precisely from the boundary
condition given by $\varphi \in Dom(\ad).$

\smallskip\noindent (6) Note that subelliptic multipliers at $x_0$
for $(p,q)$ forms are denoted by $I^q(x_0)$ without reference to
$p$, which is the holomorphic part of any such form and which
plays no role in the $\dbar$-Neumann problem.

\bigskip\noindent In the paper \cite{kohnacta}
cited above, Joseph J. Kohn proceeds by considering germs of
smooth functions at $x_0$, which he denotes by $\smooth (x_0)$. 
We will proceed now to explain Kohn's setup of his algorithm. For
this we need two more definitions:

\medskip
\newtheorem{moduledef}[subellgain]{Definition}
\begin{moduledef}
To each $x_0 \in \overline{\Omega}$ and $q \geq 1$ we associate
the module $M^q (x_0)$ defined as the set of $(1,0)$ forms
$\sigma$ satisfying that there exist a neighborhood $U$ of $x_0$
and constants $C, \epsilon > 0$ such that
\begin{equation}
||\, int(\bar\sigma)\, \varphi\, ||_{\, \epsilon}^2 \leq C \, ( \,
||\,\dbar \, \varphi\,||^2_{\, 0} + ||\, \ad  \varphi \,||^2_{\,
0} + ||\,\varphi \,||^2_{\, 0} \,) \label{moduleest}
\end{equation}
for all $(p,q)$ forms $\varphi \in \smoothc (U) \cap Dom (\ad),$
where $int(\bar\sigma)\varphi$ denotes the interior multiplication
of the $(0,1)$ form $\bar\sigma$ with the $(p,q)$ form $\varphi.$
\end{moduledef}

\smallskip\noindent The significance of $M^q (x_0)$ is that
complex gradients of subelliptic multipliers will be shown to
belong to it.

\medskip
\newtheorem{realraddef}[subellgain]{Definition}
\begin{realraddef}
Let $J \subset \smooth (x_0),$ then the real radical of $J$
denoted by $\sqrt[\R]{J}$ is the set of $g \in \smooth (x_0)$ such
that there exists some $f \in J$ and some positive natural number
$m \in \N^\ast$ such that $$|g|^m \leq |f|$$ on some neighborhood
of $x_0.$
\end{realraddef}

\smallskip\noindent The real radical is the correct generalization
of the usual radical on the ring of holomorphic functions $\oka$
for both $\C$-valued $\ra$ functions and $\C$-valued $\smooth$
functions.

We now have all definitions in place to state Kohn's Proposition
$4.7$ from \cite{kohnacta} in which he proves the properties
characterizing subelliptic multipliers that allow him to put his
algorithm together:

\medskip
\newtheorem{subellprop}[subellgain]{Proposition}
\begin{subellprop}
If $\Omega$ is a smooth pseudoconvex domain and if $x_0 \in \overline{\Omega},$
then $I^q (x_0)$ and $M^q (x_0)$ have the following properties:
\begin{enumerate}
\item[(A)] $1 \in I^n (x_0)$ and for all $q,$ whenever $x_0 \in
\Omega,$ then $1 \in I^q (x_0).$
\item[(B)] If $x_0 \in b \Omega,$ then $r \in I^q (x_0).$
\item[(C)] If $x_0 \in b \Omega,$ then $int (\theta)\, \partial
\dbar r \in M^q (x_0)$ for all smooth $(0,1)$ forms $\theta$ such
that $\langle \theta, \dbar r \rangle = 0$ on $b \Omega.$
\item[(D)] $I^q (x_0)$ is an ideal.
\item[(E)] If $f \in I^q (x_0)$ and if $g \in \smooth (x_0)$ with
$|g| \leq |f|$ in a neighborhood of $x_0,$ then $g \in I^q (x_0).$
\item[(F)] $I^q (x_0) = \sqrt[\R]{I^q (x_0)}.$
\item[(G)] $\partial I^q (x_0) \subset M^q (x_0),$ where
$\partial I^q (x_0)$ denotes the set of $(1,0)$ forms composed of
complex gradients $\partial f$ for $f \in I^q (x_0).$
\item[(H)] $\det_{n-q+1} M^q (x_0) \subset I^q (x_0),$ where $\det_{n-q+1} M^q
(x_0)$ is the coefficient of the wedge product of $n-q+1$ elements
of $M^q (x_0).$ \label{subellpropthm}
\end{enumerate}
\end{subellprop}

\smallskip\noindent {\bf Remark:} Kohn proved properties (D), (E), and (F) without employing pseudoconvexity. For properties (C) and (H), however, pseudoconvexity is a crucial hypothesis.

\bigskip\noindent Kohn gives the following corollary to Proposition~\ref{subellpropthm}, which is Theorem $1.21$ in \cite{kohnacta}:

\medskip
\newtheorem{subellcor}[subellgain]{Corollary}
\begin{subellcor}
If $\Omega$ is a smooth pseudoconvex domain and if $x_0
\in \overline{\Omega},$ then we have:
\begin{enumerate}
\item[(a)] $I^q (x_0)$ is an ideal.
\item[(b)] $I^q (x_0) = \sqrt[\R]{I^q (x_0)}.$
\item[(c)] If $r=0$ on $b \Omega,$ then $r \in I^q (x_0)$ and the
coefficients of $\partial r \wedge \dbar r \wedge (\partial \dbar
r)^{n-q}$ are in $I^q (x_0).$
\item[(d)] If $f_1, \dots, f_{n-q} \in I^q (x_0),$ then the
coefficients of $ \partial r \wedge \dbar r \wedge \partial f_1 \wedge \dots \wedge \partial f_j
\wedge (\partial \dbar
r)^{n-q-j}$ are in $I^q (x_0),$ for $j \leq
n-q.$\label{subellpropcor}
\end{enumerate}
\end{subellcor}

\smallskip\noindent {\bf Remark:} Properties (a) and (b) do not require pseudoconvexity whereas properties (c) and (d) cannot be proven in absence of pseudoconvexity. \label{pscremark}

\medskip\noindent This corollary precisely motivates how Joseph J. Kohn sets
up the algorithm.

\medskip\noindent {\bf The Kohn Algorithm:}

\medskip\noindent {\bf Step 1}  $$I^q_1(x_0) = \sqrt[\R]{(\, r,\,
\text{coeff}\{\partial r \wedge \dbar r \wedge (\partial \dbar
r)^{n-q}\}\, )}$$

\medskip\noindent {\bf Step (k+1)} $$I^q_{k+1} (x_0) = \sqrt[\R]{(\, I^q_k
(x_0),\, A^q_k (x_0)\, )},$$ where $$A^q_k (x_0)= \text{coeff}\{\partial
f_1 \wedge \dots \wedge \partial f_j \wedge \partial r \wedge
\dbar r \wedge (\partial \dbar r)^{n-q-j}\}$$ for $f_1, \dots, f_j
\in I^q (x_0)$ and $j \leq n-q.$ Note that $( \, \cdot \, )$
stands for the ideal generated by the functions inside the
parentheses and $\text{coeff}\{\partial r \wedge \dbar r \wedge (\partial
\dbar r)^{n-q}\}$ is the determinant of the Levi form for $q=1,$
namely in the $\dbar$-Neumann problem for $(0,1)$ forms.
Evidently, $I^q_k (x_0) \subset I^q (x_0)$ at each step $k,$ and
furthermore the algorithm generates an increasing chain of ideals
$$I^q_1 (x_0) \subset I^q_2 (x_0) \subset \cdots.$$ Note that the
two generators of the first ideal of multipliers $I^q_1 (x_0)$ are
globally defined objects. 

The standard setup used by Joseph J. Kohn and
others when working with the $\dbar$-Neumann problem involves a small enough neighborhood of $x_0$ that convenient frames of
vector fields and dual forms can be defined in which the boundary
condition $\varphi \in Dom(\ad)$ has a particularly simple
statement. We will adopt this type of neighborhood that Kohn
describes in section 2, page 89 of his paper \cite{kohnacta} and
describe it here for the reader's convenience. The only modification will be exchanging indices $1$ and $n$ in order to be consistent with David Catlin's setup in \cite{catlinbdry}, which we will recall starting with Section~\ref{mtbssect}.

We choose a defining function $r$ for the domain $\Omega$ such
that $|\partial r|_x = 1$ for all $x$ in a neighborhood of $b
\Omega.$ We choose a neighborhood $U$ of $x_0$ small enough that
the previous condition holds on $U,$ and we choose $(1,0)$ forms
$\omega_1, \dots, \omega_n$ on $U$ satisfying that $\omega_1 =
\partial r$ and $\langle \omega_i, \omega_j \rangle = \delta_{ij}$
for all $x \in U.$ We define by duality $(1,0)$ vector fields
$L_1, \dots, L_n$ such that $\langle \omega_i, L_j \rangle =
\delta_{ij}$ for all $x \in U.$ It follows that on $U \cap b
\Omega,$ $$L_j (r) = \bar L_j (r) = \delta_{1j}.$$ We define a
vector field $T$ on $U \cap b \Omega$ by $$T = L_1 - \bar L_1.$$
Clearly, the collection of vector fields $L_2, \dots, L_n,
\bar L_2, \dots, \bar L_n, T$ gives a local basis for the
complexified tangent space $\C T (U \cap b \Omega).$ A $(p,q)$
form $\varphi$ can be expressed in terms of the corresponding
local basis of dual forms on $U$ as $$\varphi = \sum_{|I|=p, \:
|J|=q} \: \varphi_{IJ} \, d\omega_I \wedge d \bar \omega_J,$$ for $I$
and $J$ multi-indices in $\N^n.$ As Kohn shows, $\varphi \in Dom
(\ad)$ means precisely that
$$\varphi_{IJ} (x) = 0$$ when $1 \in J$ and $x \in b \Omega.$ The
Levi form is likewise computed in this local basis. \label{prenbd}

Two of the main technical difficulties that separate the $\ra$ case solved by Kohn and the $\smooth$ case are:

\begin{enumerate}
\item[(i)] The ring $\smooth$ contains flat functions, which makes it not Noetherian.
\item[(ii)] The {\L}ojasiewicz inequalities do not hold for all elements of $\smooth (U).$
\end{enumerate}

Let us first define what it means for a
function to be flat.

\smallskip
\newtheorem{flatdef}[subellgain]{Definition}
\begin{flatdef}
A function $f \in \smooth (U)$ is said to be flat at a point $x
\in U$ if $f$ vanishes at $x$ along with its derivatives of all
orders.
\end{flatdef}

\smallskip\noindent Fortunately, the D'Angelo finite type condition guarantees that the defining function $r$ has to be non-flat in certain directions tangent to the boundary of domain $b
\Omega$ besides the non-flatness of $r$ in the normal direction guaranteed by the fact that it defines a manifold. Catlin's construction of boundary systems in \cite{catlinbdry} precisely captures this non-flatness of $r.$ He differentiates $r$ with respect to certain vector fields and their conjugates in a neighborhood of a point in such a way that he obtains real-valued functions $r_{p+2}, \dots, r_{n+1-q}$ with non-zero, linearly independent gradients, where $p$ is the rank of the Levi form at the chosen point and $b \Omega$ satisfies finite D'Angelo $q$ type. Thus, Catlin constructs non-flat elements starting from the defining function $r.$ We shall use his machinery to show this non-flatness translates into the non-flatness of the Levi determinant $\text{coeff}\{\partial r \wedge \dbar r \wedge (\partial \dbar
r)^{n-q}\}.$ Catlin's construction will be explained in detail in Section~\ref{mtbssect}.

\medskip We shall close this section by recalling from
\cite{kohnacta} the definition of the Zariski tangent space to an
ideal and to a variety, which will allow us to define the
holomorphic dimension of a variety. Catlin also defines this concept in \cite{catlinbdry}, albeit in a slightly different manner. We will also recall Catlin's definition as it is the latter we require in order to state his results from \cite{catlinbdry}.

\smallskip
\newtheorem{zariskidef}[subellgain]{Definition}
\begin{zariskidef}
Let $\I$ be an ideal in $\smooth (U)$ and let $\V (\I)$ be the
variety corresponding to $\I.$ If $x \in \V (\I),$ then we define
$Z^{\, 1,0}_x (\I)$ the Zariski tangent space of $\I$ at $x$ to be
$$Z^{\, 1,0}_x (\I) = \{ \, L \in T^{\, 1,0}_x (U) \: | \: L(f)=0 \:\:
\forall \: f \in \I \, \},$$ where $T^{\, 1,0}_x (U)$ is the
$(1,0)$ tangent space to $U \subset \C^n$ at $x.$ If $\V$ is a
variety, then we define
$$Z^{\, 1,0}_x (\V) = Z^{\, 1,0}_x (\I (\V)), $$ where $\I (\V)$
is the ideal of all functions in $\smooth (U)$ vanishing on $\V.$
\end{zariskidef}

\medskip\noindent The next lemma is Lemma $6.10$ of \cite{kohnacta}
that relates $Z^{\, 1,0}_x (\I)$ with $Z^{\, 1,0}_x (\V(\I)):$

\smallskip
\newtheorem{zariskilemma}[subellgain]{Lemma}
\begin{zariskilemma}
\label{zariski} Let $\I$ be an ideal in $\smooth (U).$ If $x \in
\V (\I),$ then
\begin{equation}
Z^{\, 1,0}_x (\V(\I)) \subset Z^{\,1,0}_x (\I). \label{zarexpr}
\end{equation}
Equality holds in ~\eqref{zarexpr} if the ideal $\I$ satisfies the
Nullstellensatz, namely $\I = \I (\V(\I)).$
\end{zariskilemma}

\medskip\noindent Let $$\nx = \{ \, L \in T^{\, 1,0}_x (b \Omega) \: |
\: \langle \, (\partial \dbar r)_x \, , \, L \wedge \bar L \,
\rangle = 0 \, \}.$$ $\nx$ is precisely the subspace of $T^{\, 1,0}_x (b \Omega)$
consisting of the directions in which the Levi form vanishes. We
end this section by defining the holomorphic dimension of a
variety sitting in the boundary of the domain $\Omega$ first as Kohn defined in \cite{kohnacta} and then as Catlin defined it in \cite{catlinbdry}:

\smallskip
\newtheorem{holdimdefk}[subellgain]{Definition (Kohn)}
\begin{holdimdefk}
Let $\V$ be a variety in $U$ corresponding to an ideal $\I$ in
$\smooth (U)$ such that $\V \subset b \Omega.$ We define the
holomorphic dimension of $\V$ in the sense of Kohn by $$hol. \: dim \, (\V) = \min_{x
\in \V} \dim Z^{\, 1,0}_x (\V) \cap \nx.$$ \label{holdimkohn}
\end{holdimdefk}

\smallskip
\newtheorem{holdimdefc}[subellgain]{Definition (Catlin)}
\begin{holdimdefc}
Let $\V$ be a variety in $U$ corresponding to an ideal $\I$ in
$\smooth (U)$ such that $\V \subset b \Omega.$ We define the
holomorphic dimension of $\V$ in the sense of Catlin by $$hol. \: dim \, (\V) = \max_{x
\in \V} \dim Z^{\, 1,0}_x (\V) \cap \nx.$$ \label{holdimcatlin}
\end{holdimdefc}

\section{Finite D'Angelo type}
\label{findap}

Starting with \cite{opendangelo}, John D'Angelo introduced various numerical functions that measure the maximum order of contact of holomorphic varieties of complex dimension $q$ with a real hypersurface in $\C^n$ such as the boundary of a domain. The interested reader should consult \cite{dangelo} for the most comprehensive treatment of this topic.

We shall first give the classical definition of order of contact for $q=1,$ holomorphic curves, where the most natural definition is clear. We shall then discuss ways to understand this order of contact for $q>1.$ Let $\cv=\cv(m,p)$ be the set of all germs of holomorphic curves $$\varphi: (U,0) \rightarrow (\C^m, p),$$ where $U$ is some neighborhood of the origin in $\C^1$ and $\varphi(0)=p.$ For all $t \in U,$ $\varphi(t)=(\varphi_1(t), \dots, \varphi_m(t)),$ where $\varphi_j(t)$ is holomorphic for every $j$ with $1 \leq j \leq m.$ For each component $\varphi_j,$ the order of vanishing at the origin ${\text ord}_0 \, \varphi_j$ is the order of the first non-vanishing derivative of $\varphi_j,$ i.e. $s \in \N$ such that $$\frac{d}{dt} \varphi_j (0) = \cdots = \frac{d ^{s-1}}{dt^{s-1}} \varphi_j (0) = 0,$$ but $\frac{d ^s}{dt^s} \varphi_j (0) \neq 0.$ We set ${\text ord}_0 \, \varphi = \min_{1 \leq j \leq m} \, {\text ord}_0 \, \varphi_j.$ Consider $\varphi^* r,$ the pullback of $r$ to $\varphi,$ and let $ {\text ord}_0 \, \varphi^* r$ be the order of the first non-vanishing derivative at the origin of $ \varphi^* r$ viewed as a function of $t.$

\medskip
\newtheorem{firstft}{Definition}[section]
\begin{firstft}
Let $M$ be a real hypersurface in $\C^n,$ and let $r$ be a defining function for $M.$ The D'Angelo $1$-type at $x_0 \in M$ is given by $$\Delta_1 (M, x_0) = \sup_{\varphi \in \cv(n,x_0)} \frac{ {\text ord}_0 \, \varphi^* r}{{\text ord}_0 \, \varphi}.$$If $\Delta_1 (M, x_0)$ is finite, we call $x_0$ a point of finite D'Angelo $1$-type.
\end{firstft}

\medskip\noindent {\bf Remark:} John D'Angelo showed $\Delta_1 (M, x_0)$ is well-defined, i.e. independent of the defining function $r$ chosen for $M.$ He also showed it fails to be upper semi-continuous in \cite{nonusc}. Fortunately,  $\Delta_1 (M, x_0)$ being finite is an open condition, and  $\Delta_1 (M, x_0)$ itself is finitely determined, i.e. it is determined merely by a certain jet of the Taylor expansion of $r$ at $x_0$ and not the entire Taylor expansion. We shall rigorously state both of these properties after we define  $\Delta_q (M, x_0).$

\smallskip When holomorphic varieties have complex dimension greater than $1,$ there is no longer just one natural definition of their order of contact with a real hypersurface in $\C^n.$ One approach is to reduce this case to computing $\Delta_1 (\tilde M, x_0)$ for a related hypersurface $\tilde M$ sitting in a different $\C^m.$ This is the approach taken by D'Angelo in \cite{opendangelo}. Let $\phi : \C^{n-q+1} \rightarrow \C^n$ be any linear embedding of $\C^{n-q+1}$ into $\C^n.$ For generic choices of $\phi,$ the pullback $\phi^* M$ will be a hypersurface in $\C^{n-q+1}.$ We can thus define $\Delta_q (M, x_0)$ as follows:

\medskip
\newtheorem{qft}[firstft]{Definition}
\begin{qft}
Let $M$ be a real hypersurface in $\C^n,$ and let $r$ be a defining function for $M.$ The D'Angelo $q$-type at $x_0 \in M$ is given by $$\Delta_q (M, x_0) =\inf_\phi \sup_{\varphi \in \cv(n-q+1,x_0)} \frac{ {\text ord}_0 \, \varphi^* \phi^*r}{{\text ord}_0 \, \varphi}=\inf_\phi \Delta_1 (\phi^*r, x_0),$$ where $\phi : \C^{n-q+1} \rightarrow \C^n$ is any linear embedding of $\C^{n-q+1}$ into $\C^n$ and we have identified $x_0$ with $\phi^{-1} (x_0).$ If $\Delta_q (M, x_0)$ is finite, we call $x_0$ a point of finite D'Angelo $q$-type.
\end{qft}

By truncating the Taylor expansion of the defining function $r$ at $x_0,$ John D'Angelo was able to prove one of the most important properties of $\Delta_q (M, x_0),$ namely the openness of the set of points of finite $q$-type. Furthermore, using ideas from algebraic geometry over the ring of holomorphic functions, John D'Angelo was able to give a bound on the maximal jump of $\Delta_q (M, x)$ in a neighborhood of $x_0.$ We recall Theorem 6.2 from p.634 of \cite{opendangelo}:

\medskip
\newtheorem{openqt}[firstft]{Theorem}
\begin{openqt}
Let $M$ be a smooth real hypersurface in $\C^n$ and let $\Delta_q (M, x_0)$ be finite at some $x_0 \in M,$ then there exists a neighborhood $V$ of $x_0$ on which $$\Delta_q (M, x) \leq 2 (\Delta_q(M, x_0))^{n-q}.$$ \label{openqtthm}
\end{openqt}

\medskip\noindent {\bf Remark:} This theorem holds {\bf independently} of pseudoconvexity, which will be crucial for the proof of the Main Theorem~\ref{maintheorem}.

\medskip\noindent We shall now state another essential property of $\Delta_q (M, x_0),$ being finite determined. This is Proposition 14 from p.88 of \cite{dangeloftc}:

\medskip
\newtheorem{findetqt}[firstft]{Theorem}
\begin{findetqt}
The function $\Delta_q(M, x_0)$ is finite determined. In other words, if $\Delta_q(M, x_0)$ is finite, then there exists an integer $k$ such that $\Delta_q(M, x_0)=\Delta_q(M',x_0)$ for $M'$ a hypersurface defined by any $r'$ that has the same $k$ jet at $x_0$ as the defining function $r$ of $M.$ \label{findetqtthm}
\end{findetqt}

\medskip\noindent {\bf Remark:} Let $t = \Delta_q (M, x_0) < \infty,$ then it follows from the proof of Proposition 14 in \cite{dangeloftc} that we can let $k = \lceil t \rceil,$ the roundup of $t,$ i.e. the lowest integer greater than or equal to $t.$

\medskip\noindent We shall close the section with a brief discussion of another approach to defining the order of contact of holomorphic varieties that have complex dimension greater than $1$ taken by David Catlin in \cite{catlinsubell}.

David Catlin wished to avoid having to characterize the order of contact of a holomorphic variety $V^q$ of complex dimension $q$ with the boundary of the domain along the singular locus of the variety, which can be considerably more complicated when $q>1$ than for holomorphic curves. To that end, he introduced a numerical function $D_q (M, x_0)$ that measures the order of contact of varieties $V^q$ with $M$ only along generic directions.

Let $V^q$ be the germ of a holomorphic variety of complex dimension $q$ passing through $x_0.$ Let $W$ be the set of all $(n-q+1)$-dimensional complex planes through $x_0.$ Consider the intersection $V^q \cap S$ for $S \in W.$ For a generic, thus open and dense, subset $\tilde W$ of $W,$ $V^q \cap S$ consists of finitely many one-dimensional components $V^q_{S,k}$ for $k=1, \dots, P.$ Let us parametrize these curves by some open set $U \ni 0$ in $\C.$ Thus, $\gamma_S^k : U \rightarrow V^q_{S,k},$ where $\gamma_S^k (0)=x_0.$ Set $$\tau (V^q \cap S, x_0) = \max_{k=1, \dots, P} \frac{ {\text ord}_0 \, {\left(\gamma^k_S\right)}^* r}{{\text ord}_0 \, \gamma^k_S}.$$ In Section 3 of \cite{catlinsubell}, David Catlin showed $\tau (V^q \cap S, x_0)$ assumes the same value for all $S$ in a generic subset $\tilde W$ of planes. Thus he defined $$\tau(V^q, x_0) = {\text gen.val} \left\{ \tau (V^q \cap S, x_0) \right\}.$$

\medskip
\newtheorem{catlinft}[firstft]{Definition}
\begin{catlinft}
Let $M$ be a real hypersurface in $\C^n.$ The Catlin $q$-type at $x_0 \in M$ is given by $$D_q (M, x_0)=\sup_{V^q}  \left\{ \tau (V^q, x_0) \right\},$$ where the supremum is taken over the set of all germs of $q$-dimensional holomorphic varieties $V^q$ passing through $x_0.$
\end{catlinft}

\noindent Clearly, $\Delta_1 (M, x_0) = D_1 (M , x_0),$ but for $q>1,$ the relationship between $\Delta_q (M, x_0)$ and $D_q (M, x_0)$ is not so clear. It is known, however, that $\Delta_q (M, x_0)$ is finite iff $D_q (M, x_0)$ is finite. The interested reader should consult the beginning of Section 6 of \cite{kohndangelo}.

\section{Catlin's multitype and boundary systems}
\label{mtbssect}

This section is devoted to recalling the concepts of boundary system, multitype, and commutator multitype from David Catlin's paper \cite{catlinbdry}. Let $x_0 \in b \Omega.$ A boundary system $$\B_\nu = \{r_1, r_{p+2}, \dots, r_\nu; L_2, \dots, L_\nu\}$$ of rank $p$ and codimension $n-\nu$ is a collection of $\R$-valued smooth functions with linearly independent gradients and $(1,0)$ vector fields in a neighborhood of $x_0$ in $\C^n.$ The Levi form has rank $p$ at $x_0,$ $r_1=r$ is the defining function of the domain, and the other functions $ r_{p+2}, \dots, r_\nu$ are obtained from $r$ by differentiation with respect to the vector fields $ L_2, \dots, L_\nu$ and their conjugates in a certain order that will be described in detail shortly. The commutator type  $\ct(x_0)$ is an $n$-tuple of positive rational numbers or $+\infty$ whose first $\nu$ entries measure how many times $r$ has to be differentiated via $L_2, \dots, L_\nu$ and their conjugates until the functions  $ r_{p+2}, \dots, r_\nu$ are obtained. The multitype $\mt(x_0)$ is likewise an  $n$-tuple of positive rational numbers or $+\infty$ that measures the non-vanishing of the defining function $r$ in different directions. The multitype  $\mt(x_0)$ and the commutator multitype $\ct(x_0)$ equal each other when the domain is pseudoconvex. It is not known known whether this equality holds more generally. It should be noted that for a pseudoconvex domain of finite D'Angelo $q$-type at $x_0,$  $\B_\nu$ exists for $\nu=n+1-q$ and both the multitype $\mt(x_0)$ and the commutator multitype $\ct(x_0)$ have only finite entries up to their $\nu^{th}$ one.

Among the concepts of boundary system, multitype, and commutator multitype, the easiest notion to introduce is that of multitype, so we will follow Catlin in \cite{catlinbdry} in describing it first. Since both $\mt(x_0)$ and $\ct(x_0)$ are $n$-tuples of rational numbers or $+\infty$, we need to consider all such $n$-tuples that satisfy certain properties. We will call these weights. We will first define these weights and then describe a subset of weights with even better properties, which we will call the set of distinguished weights.

\medskip
\newtheorem{wgt}{Definition}[section]
\begin{wgt}
\label{wgtdef} Let $\Gamma_n$ denote the set of $n$-tuples of rational numbers $\Lambda=(\lambda_1, \dots, \lambda_n)$ with $1 \leq \lambda_i \leq + \infty$ satisfying the following two properties:
\begin{enumerate}
\item[(i)] $\lambda_1 \leq \lambda_2 \leq \cdots \leq \lambda_n.$
\item[(ii)] For each $k$ such that $1 \leq k \leq n,$ either $\lambda_k = + \infty$ or there exists a set of integers $a_1, \dots , a_k$ such that $a_j >0$ for all $1 \leq j \leq k$ and $$\sum_{j=1}^k \frac{a_j}{\lambda_j}=1.$$
\end{enumerate}
The set $\Gamma_n$ is ordered lexicographically, i.e. given $\Lambda', \Lambda'' \in \Gamma_n$ such that $\Lambda'=(\lambda'_1, \dots, \lambda'_n)$ and $\Lambda''=(\lambda''_1, \dots, \lambda''_n),$ then $\Lambda' < \Lambda''$ if there exists $k$ with $1 \leq k \leq n$ such that $\lambda'_j = \lambda''_j$ for all $j<k$ and $\lambda'_k < \lambda''_k.$ The set $\Gamma_n$ is called the set of weights.
\end{wgt}

\smallskip\noindent {\bf Remark:}
Requiring the sum to equal $1$ at each step $k$ in property (ii) is one of Catlin's most remarkable ideas as it enables him to prove the equality of the multitype and the commutator type for a pseudoconvex domain by truncating the defining function $r$ with respect to a weight $\Lambda=(\lambda_1, \dots, \lambda_n)$ in a way that preserves all terms in the Taylor expansion of $r$ that are essential for producing the functions $r_{p+2}, \dots, r_\nu$ under differentiation. 

\medskip\noindent Let us now define distinguished weights and the multitype $\mt(x_0):$

\medskip
\newtheorem{diswgt}[wgt]{Definition}
\begin{diswgt}
\label{diswgtdef} Let $\Omega \subset \C^n$ be a smooth domain with defining function $r.$ A weight $\Lambda=(\lambda_1, \dots, \lambda_n) \in \Gamma_n$ is called distinguished if there exist holomorphic coordinates $(z_1, \dots, z_n)$ around $x_0$ such that 
\begin{enumerate}
\item[(i)] $x_0$ is mapped at the origin;
\item[(ii)] If $\sum_{i=1}^n \: \frac{\alpha_i+ \beta_i}{\lambda_i}<1,$ then $D^\alpha \bar D^\beta r(0)=0,$ where $D^\alpha= \frac{\partial^{|\alpha|}}{\partial z^{\alpha_1}_1 \cdots \partial z^{\alpha_n}_n}$ and $\bar D^\beta= \frac{\partial^{|\beta|}}{\partial \bar z^{\beta_1}_1 \cdots \partial \bar z^{\beta_n}_n}.$
\end{enumerate}
We will denote by $\tilde \Gamma_n (x_0)$ the set of distinguished weights at $x_0.$
\end{diswgt}

\smallskip\noindent {\bf Remark:} Property (ii) in Definition~\ref{wgtdef} and property (ii) in Definition~\ref{diswgtdef} taken together show that the underlying idea of this setup is to measure the order of vanishing of the defining function $r$ in various directions. It turns out this measure is a weight in $\Gamma_n$ called the multitype. We clearly have to measure the vanishing of $r$ in a way that is independent of the local coordinates chosen. This precisely justifies the wording of the next definition.

\medskip
\newtheorem{mlt}[wgt]{Definition}
\begin{mlt}
\label{mltdef} The multitype $\mt(x_0)$ is defined to be the smallest weight in lexicographic sense $\mt(x_0) = (m_1, \dots, m_n)$ such that $\mt(x_0) \geq \Lambda$ for every distinguished weight $\Lambda \in \tilde \Gamma_n (x_0).$
\end{mlt}

\medskip\noindent We will now state the main theorem of Catlin's paper \cite{catlinbdry} from page 531 that summarizes the properties of the multitype $\mt(x_0):$

\medskip
\newtheorem{mltthm}[wgt]{Theorem}
\begin{mltthm}
\label{mlttheorem} Let $\Omega \subset \C^n$ be a pseudoconvex domain with smooth boundary. Let $x_0 \in b \Omega.$ The multitype $\mt(x_0)$ has the following properties:
\begin{enumerate}
\item[(1)] $\mt(x_0)$ is upper semi-continuous with respect to the lexicographic ordering, i.e. there exists a neighborhood $U \ni x_0$ such that for all $x \in U \cap b \Omega,$ $\mt(x) \leq \mt(x_0).$
\item[(2)] If $\mt(x_0)=(m_1, \dots, m_n)$ satisfies that $m_{n-q}<\infty,$ then there exist a neighborhood $U \ni x_0$ and a submanifold $M$ of $U \cap b \Omega$ of holomorphic dimension at most $q$ in the sense of Catlin such that $x_0 \in M$ and the level set of $\mt(x_0)$ satisfies $$\{ x \in U \cap b \Omega \: \big| \: \mt(x)=\mt(x_0)\} \subset M.$$
\item[(3)] If $\mt(x_0)=(m_1, \dots, m_n),$ then there exist coordinates $(z_1, \dots, z_n)$ around $x_0$ such that $x_0$ is mapped to the origin and if $\sum_{i=1}^n \: \frac{\alpha_i+ \beta_i}{m_i}<1,$ then $D^\alpha \bar D^\beta r(0)=0.$ If one of the entries $m_i = + \infty$ for some $1 \leq i \leq n,$ then these coordinates should be interpreted in the sense of formal power series.
\item[(4)] If $\mt(x_0)=(m_1, \dots, m_n),$ then for each $q=1, \dots, n,$ $$m_{n+1-q} \leq \Delta_q (b \Omega, x_0),$$ where $\Delta_q (b \Omega, x_0)$ is the D'Angelo $q$ type of the point $x_0,$ i.e. the maximum order of contact of varieties of complex dimension $q$ with the boundary of $\Omega$ at $x_0.$
\end{enumerate}
\end{mltthm}

\medskip\noindent It is clear that Definition~\ref{mltdef} does not specify a procedure for computing $\mt(x_0)$ for a domain $\Omega$ at the boundary point $x_0.$ Instead, Catlin defined another weight $\ct(x_0) \in \Gamma_n$ called the commutator multitype, which he proceeded to compute by differentiating $r.$ In the process of computing $\ct(x_0),$ Catlin came up with the definition of a boundary system. He then showed that $\ct(x_0)=\mt(x_0)$ for a pseudoconvex domain. We will now explain his construction of the commutator multitype $\ct(x_0)$ and of a boundary system $\B_\nu (x_0)= \{r_1, r_{p+2}, \dots, r_\nu; L_2, \dots, L_\nu\}.$

The commutator multitype  $\ct(x_0)=(c_1, \dots, c_n) \in \Gamma_n$ always satisfies that $c_1=1$ because as explained on page~\pageref{prenbd}, $L_1(r)=1,$ which comes from the fact that $r$ describes a manifold and thus its gradient has a non-zero component in the normal direction. This means only one differentiation of $r$ in the direction of $L_1$ suffices to produce a non-vanishing function at $x_0,$ hence $c_1=1.$ We set $r_1=r.$ Next, suppose that the Levi form of $b \Omega$ at $x_0$ has rank equal to $p.$ In this case, set $c_i=2$ for $i=2, \dots, p+1.$ Without loss of generality, we can choose in the construction on page~\pageref{prenbd} the smooth vector fields of type $(1,0)$ $L_2, \dots, L_{p+1}$ such that $L_i(r)=\partial r (L_i)\equiv 0$ and the $p \times p$ Hermitian matrix $\partial \bar \partial r (L_i, L_j)(x_0)$ for $2 \leq i,j \leq p+1$ is nonsingular. The reader should note that round parentheses stand for the evaluation of forms on vector fields. If $p+1 \geq \nu,$ we have finished the construction of the boundary system $\B_\nu(x_0).$

If $p+1 < \nu,$ we need to explain next how the rest of the vector fields $L_{p+2}, \dots, L_\nu$ and the functions $r_{p+2}, \dots, r_\nu$ are chosen in the boundary system $\B_\nu (x_0).$ Let us consider the $(1,0)$ smooth vector fields in the kernel of the Levi form at $x_0.$ We thus denote by $T^{\, (1,0)}_{p+2}$ the bundle consisting of $(1,0)$ vector fields $L$ such that $\partial r (L)=0$ and $\partial \bar \partial r (L, \bar L_j)=0$ for $j=2, \dots, p+1.$ We follow Catlin in passing to the set of germs of sections of $T^{\, (1,0)}_{p+2},$ which we will denote by $\T_{p+2}.$ Germs are not technically necessary for defining a boundary system and the commutator multitype, but they become essential \label{germsgood} later on when Catlin describes truncated boundary systems because one can then pick a representative in the equivalence class of a germ given by a vector field with polynomial coefficients.

All the directions in which the defining function vanishes up to order $2$ have already been identified. It is thus clear we have to consider next lists of vector fields of length at least $3.$ Let $l \in \N$ be such that $l \geq 3.$ Denote by $\ls$ a list of vector fields $\ls= \{L^1, \dots, L^l\}$ such that there is a fixed, non-vanishing vector field $L \in T^{\, (1,0)}_{p+2}$ and $L^i = L$ or $L^i= \bar L$ for all $1 \leq i \leq l.$ Let $\ls  \partial r$ be the function $$\ls  \partial r (x) = L^1 \cdots L^{l-2} \, \partial r \, ([L^{l-1},L^l]) (x)$$ for $x \in b \Omega.$ The reader should note that if both $L^{l-1}$ and $L^l$ are $L$ or both of them are $\bar L,$ then the commutator $[L^{l-1},L^l]$ vanishes. Thus, to have any chance of obtaining a non-vanishing function as $\ls \partial r,$ one of $L^{l-1}$ and $L^l$ should be $L$ and the other one $\bar L.$ Therefore, $\ls \partial r (x_0)$ measures the vanishing order of the diagonal entry of the Levi form at $x_0$ corresponding to $L.$ We distinguish two cases:

\smallskip\noindent {\bf Case 1:} If $\ls \partial r (x_0)=0$ for every such list $\ls,$ then set $c_{p+2}= \infty.$ Given that weights are increasing $n$-tuples of rational numbers, it follows $c_i = \infty$ for all $i=p+2, \dots, n.$ We have finished the construction of both the boundary system $\B_\nu(x_0)$ and of the commutator multitype $\ct (x_0).$

\smallskip\noindent {\bf Case 2:} There exists at least one list $\ls$ such that $\ls \partial r (x_0)\neq 0.$ Among all lists with this property, we choose one list for which the length $l$ is the smallest. Clearly, there might exist more than one list of smallest length, but the entries $c_i$ of the commutator multitype $\ct(x_0)$ will turn out to be independent of the choice made here. Set $c_{p+2} = l,$ where $l$ is this smallest value of the length of the list. Let $\ls_{p+2} = \{ L^1, \dots, L^l \}$ be the list chosen whose length satisfies $l = c_{p+2}$ and let $L \in \T_{p+2}$ be the germ of the fixed vector field in $T^{\, (1,0)}_{p+2}$ such that $L^i=L$ or $L^i= \bar L$ for all $1 \leq i \leq l.$ Define functions $f$ and $g$ by $$f (x) =Re \{ L^2 \cdots L^{l-2} \, \partial r \, ([L^{l-1},L^l]) (x)\}$$ and $$g (x) =Im \{ L^2 \cdots L^{l-2} \, \partial r \, ([L^{l-1},L^l]) (x)\}.$$ Since $l \geq 3,$ the definitions of $f$ and $g$ make sense. We define $\R$-valued vector fields $X$ and $Y$ such that $L=X+iY.$ Since $l$ was chosen to be minimal, it follows that $f(x_0)=g(x_0)=0,$ but $L^1(f+ig)(x_0) \neq 0.$ This implies at least one of $Xf(x_0),$ $Xg(x_0),$ $Yf(x_0),$ and $Yg(x_0)$ does not vanish. Without loss of generality, let $Xf(x_0) \neq 0.$ We set $r_{p+2}(x)=f(x)$ and $L_{p+2}=L,$ the vector field used in constructing the list $\ls_{p+2}.$ It follows that $L_{p+2}(r_{p+2})(x_0) \neq 0.$ This concludes the second case and thus the construction at step $p+2.$

We proceed inductively. Assume that for some integer $\nu-1$ with $p+2 \leq \nu-1<n$ we have already constructed finite positive numbers $c_1, \dots, c_{\nu-1}$ as well as functions $r_1,$ $r_{p+2}, \dots, r_{\nu-1}$ and vector fields $L_2, \dots, L_{\nu-1}.$ Let $T^{\, (1,0)}_\nu$ denote the set of $(1,0)$ smooth vector fields $L$ such that $\partial \bar \partial r (L, \bar L_j)=0$ for $j=2, \dots, p+1$ and $L(r_k)=0$ for $k=1,p+2, p+3, \dots,\nu-1.$ \label{vfbehavior}  Let $\T_\nu$ be the set of germs of sections of $T^{\, (1,0)}_\nu.$ For each $k=p+2, \dots,\nu-1,$ $L_k \in \T_k$ and $L_k(r_k)(x_0) \neq 0,$ which implies that $T^{\, (1,0)}_\nu$ is a subbundle of $T^{\, (1,0)} (b \Omega)$ of dimension $n+1-\nu$ because the vector fields $L_2, \dots, L_{\nu-1}$ were chosen to be linearly independent. We have to describe next the list $\ls$ for which we will compute $\ls \partial r (x).$ We are allowed to use both vector fields from among $L_{p+2}, \dots, L_{\nu-1}$ as well as vector fields in $\T_\nu.$ Thus, we fix some vector field $L$ in $\T_\nu$ and consider the list $\ls=\{L^1, \dots, L^l\}$ such that each $L^i$ is one of the vector fields from the set $\{L_{p+2}, \bar L_{p+2}, \dots, L_{\nu-1}, \bar L_{\nu-1},L, \bar L\}.$ Let $l_i$ denote the total number of times both $L_i$ and $\bar L_i$ occur in $\ls$ for $p+2 \leq i \leq \nu-1$ and let $l_\nu$ denote the total number of times both $L$ and $\bar L$ occur in the list $\ls.$ We now introduce two definitions that pertain to the list $\ls$ and explain their significance:

\medskip
\newtheorem{ordlist}[wgt]{Definition}
\begin{ordlist}
\label{ordlistdef} A list $\ls= \{L^1, \dots, L^l\}$ is called ordered if
\begin{enumerate}
\item[(i)] $L^j=L$ or $L^j= \bar L$ for $1 \leq j \leq l_\nu$
\item[(ii)] $L^j=L_i$ or $L^j=\bar L_i$ for $1+\sum_{k=i+1}^\nu l_k \leq j \leq \sum_{k=i}^\nu l_k.$
\end{enumerate}
\end{ordlist}

\smallskip\noindent {\bf Remarks:} 
\begin{enumerate}
\item[(1)] Part (i) says that the differentiation with respect to the extra vector field $L \in \T_\nu$ or its conjugate should be done outside of any differentiation with respect to the previously chosen vector fields $L_{p+2}, \dots, L_{\nu-1}.$ 
\item[(2)] Part (ii) of the definition says that if we look from left to right at the list $\ls,$ we should have differentiation with respect to $L$ or $\bar L,$ then differentiation with respect to $L_{\nu-1}$ or its conjugate if it takes place at all, then differentiation with respect to $L_{\nu-2}$ or its conjugate if it takes place, and so on. Therefore, we differentiate $r$ with respect to previously chosen vector fields inside and with respect to the more recently chosen vector fields outside. Not allowing vector fields to mix makes vanishing orders of $r$ easier to understand and exploit.
\end{enumerate}

\medskip
\newtheorem{admlist}[wgt]{Definition}
\begin{admlist}
\label{admlistdef} A list $\ls= \{L^1, \dots, L^l\}$ is called $\nu$-admissible if
\begin{enumerate}
\item[(i)] $l_\nu>0;$
\item[(ii)] $$\sum_{i=p+2}^{\nu-1} \frac{l_i}{c_i}<1,$$ where $\ct^{\nu-1}=(c_1, \dots, c_{\nu-1})$ is the $(\nu-1)^{th}$ commutator multitype.
\end{enumerate}
\end{admlist}

\smallskip\noindent {\bf Remarks:} 
\begin{enumerate}
\item[(1)] It is obvious that condition (i) should be imposed because if we do not differentiate with respect to the new vector field $L \in \T_\nu$ or its conjugate, then we cannot expect to obtain anything that we did not already have by step $\nu-1.$
\item[(2)] Condition (ii) follows from the minimality of the length of the lists chosen at the previous steps. In other words, if we strip away the first $l_\nu$ vector fields from $\ls$ and look at $\ls'=\{L^{l_\nu+1}, \dots, L^l\},$ then this is a list that appeared at one of the previous steps, so we know $\ls' \partial r(x_0)=0$ if (ii) holds.
\item[(3)] It also makes perfect sense that lists with property (ii) should be considered since we are trying to construct a commutator multitype $\ct(x_0)$ such that $\ct(x_0)=\mt(x_0),$ and the multitype $\mt(x_0)$ is defined as the weight that dominates all distinguished weights.
\end{enumerate}

\smallskip\noindent We now consider only $\nu$-admissible, ordered lists $\ls$ and distinguish two cases:

\smallskip\noindent {\bf Case 1:} For all such lists $\ls,$ $\ls \partial r_1 (x_0)=0.$ In this case, we set $c_\nu=\infty.$ It follows that $c_i= \infty$ for all $\nu \leq i \leq n.$ We have finished the construction of both the boundary system $\B_\nu(x_0)$ and of the commutator multitype $\ct (x_0).$

\smallskip\noindent {\bf Case 2:} There exists at least one such list $\ls$ for which $\ls \partial r_1 (x_0) \neq 0.$ We would like to choose the list $\ls$ with minimal length just as before. Let $c(\ls)$ denote the solution to the equation $$\sum_{i=p+2}^{\nu-1} \frac{l_i}{c_i}+\frac{l_\nu}{c(\ls)}=1.$$ Because $\ls$ is $\nu$-admissible and thus satisfies condition (ii) of Definition~\ref{admlistdef}, $$1-\sum_{i=p+2}^{\nu-1} \frac{l_i}{c_i}>0,$$ and $1-\sum_{i=p+2}^{\nu-1} \frac{l_i}{c_i}$ is a rational number since all entries $c_i$ are rational and all numbers $l_i$ are positive integers. $l_\nu$ is also a positive integer, so the solution $c(\ls)$ has to be a positive rational number. Set $$c_\nu=\inf \{c(\ls) \: \big| \: \ls \: \text{is} \: \nu-\text{admissible, ordered, and satisfies} \: \ls \partial r_1 (x_0) \neq 0 \}.$$ If there exists more than one such list for which $c(\ls)$ reaches the infimum, we make an arbitrary choice and denote it by $\ls_\nu=\{L^1, \dots, L^l\}.$ Next, we set $\ls'_\nu=\{L^2, \dots, L^l\},$ and define $\R$-valued vector fields $X$ and $Y$ such that $L^1=X+iY.$ Just as before, we let functions $f$ and $g$ be defined by $f(x)=Re\{\ls'_\nu \partial r_1(x)\}$ and $g(x)=Im\{\ls'_\nu \partial r_1(x)\}$ The minimality of the length of $\ls_\nu$ along with part (i) of Definition~\ref{ordlistdef} and part (ii) of Definition~\ref{admlistdef} together imply that $f(x_0)=g(x_0)=0.$ Since $\ls_\nu \partial r_1(x_0)= (X+iY)(f+ig)(x_0)\neq 0,$ it follows that at least one of $Xf(x_0),$ $Xg(x_0),$ $Yf(x_0),$ and $Yg(x_0)$ does not vanish. Without loss of generality, let $Xf(x_0) \neq 0.$ We set $r_\nu(x)=f(x)$ and $L_\nu=L,$ the vector field in $\T_\nu$ used in composing the list $\ls_\nu.$ This concludes the second case as well as the construction of the boundary system at step $\nu.$

\medskip\noindent {\bf Remarks:}
\begin{enumerate}
\item[(1)] If $\ls=\{L^1, \dots, L^l\}$ is any ordered list such that $L^i=L_j$ or $L^i=\bar L_j$ for all $1 \leq i \leq l$ and all $p+2 \leq j \leq \nu$ and if $\sum_{i=p+2}^{\nu} \frac{l_i}{c_i}<1,$ then $\ls \partial r_1 (x_0)=0.$ In the case $l_\nu >0,$ this follows from the minimality of $c_\nu.$ In the case $l_\nu =0,$ this follows from the minimality of the previously chosen $c_{p+2}, \dots, c_{\nu-1}.$
\item[(2)] The construction of the boundary system $\B_\nu (x_0) = \{r_1, r_{p+2}, \dots, r_\nu; L_2, \dots, L_\nu\}$ depends on the choices of lists $\ls_{p+2}, \dots, \ls_\nu$ as well as of vector fields $L_{p+2}, \dots, L_\nu$ in these lists. Fortunately, by minimality of the lengths of the lists chosen, the commutator multitype $\ct(x_0)=(c_1, \dots, c_n)$ is the same regardless of these choices. 
\end{enumerate}

\medskip\noindent We shall call a collection $$\B_\nu (x_0) = \{r_1, r_{p+2}, \dots, r_\nu; L_2, \dots, L_\nu\}$$ of functions and vector fields a boundary system of rank $p$ and codimension $n-\nu$ if it is obtained by the procedure described above. The vector fields $L_2, \dots, L_\nu$ are called the special vector fields associated to the boundary system $\B_\nu.$ The reader should note that while $\ct(x_0)$ and $\mt(x_0)$ always exist, $\B_\nu$ only exists for $\nu \geq 2$ if there are $\nu-1$ vector fields $L_2,  \dots, L_\nu$ that can be chosen according to the procedure outlined above. 

The $\nu^{th}$ commutator multitype of the boundary system $\B_\nu$ is the $\nu$-tuple $\ct^\nu=(c_1, \dots, c_\nu).$ We summarize in the next theorem two of the most important properties of $\ct^\nu,$ which are contained in Proposition 2.1 on page 536 and Theorem 2.2 on page 538 of Catlin's paper \cite{catlinbdry}:

\medskip
\newtheorem{cmtype}[wgt]{Theorem}
\begin{cmtype}
\label{cmtypethm} Let $\Omega= \{ z \in \C^n \: \big| \: r(z)<0\}$ be a smoothly bounded domain, and let $x_0 \in b \Omega.$ The $\nu^{th}$ commutator multitype $\ct^\nu=(c_1, \dots, c_\nu)$ of the boundary system $\B_\nu$ satisfies the following two properties:
\begin{enumerate}
\item[(i)] $\ct^\nu$ is upper semi-continuous with respect to the lexicographic ordering, i.e. there exists a neighborhood $U \ni x_0$ such that for all $x \in U \cap b \Omega,$ $\ct^\nu (x) \leq \ct^\nu (x_0).$
\item[(ii)] If $\Omega$ is pseudoconvex, then $\ct(x_0) = \mt(x_0),$ so $\ct^\nu(x_0)=\mt^\nu(x_0),$ where $\mt^\nu=(m_1, \dots, m_\nu)$ consists of the first $\nu$ entries of the multitype $\mt=(m_1, \dots, m_n).$
\end{enumerate}
\end{cmtype}

\smallskip\noindent Next, we would like to understand the stratification induced by the partial commutator multitype $\ct^\nu.$  We start this discussion by stating Proposition 2.1 on page 536 of Catlin's paper \cite{catlinbdry} strengthened in an obvious manner. The differences between this statement and Catlin's original statement will be outlined in a remark following the proposition.

\medskip
\newtheorem{stratcatlin}[wgt]{Proposition}
\begin{stratcatlin}
\label{stratcatlinprop} Let $\B_\nu$ for $p+2 \leq \nu \leq n$ be a boundary system of rank $p$ and codimension $n- \nu$ in a neighborhood of a given boundary point $x_0.$ There exists a neighborhood $U$ of $x_0$ such that all the following conditions are satisfied on its closure $\overline{U}$:
\begin{enumerate}
\item[(i)] For all $x \in \overline{U} \cap b \Omega,$ $\ct^\nu (x) \leq \ct^\nu (x_0),$ where $\ct^\nu=(c_1, \dots, c_\nu)$ is the $\nu^{th}$ commutator multitype;
\item[(ii)] $$M^\nu=\{ x \in \overline{U} \cap b \Omega \: \big| \: r_j(x)=0, \: j=1,p+2, \dots, \nu\}$$ is a submanifold of $\overline{U} \cap b \Omega$ of holomorphic dimension $n-\nu$ in the sense of Catlin;
\item[(iii)] The level set of the commutator multitype at $x_0$ satisfies that $$\{ x \in \overline{U} \cap b \Omega \: \big| \: \ct^\nu(x)=\ct^\nu(x_0)\} \subset M^\nu;$$
\item[(iv)] For all $x \in \overline{U} \cap b \Omega,$ the Levi form has rank at least $p$ at $x;$
\item[(v)] For all $x \in \overline{U} \cap b \Omega,$ $\ls_j \partial r_1(x) \neq 0$ for all $j=p+2, \dots, \nu,$ where $\ls_{p+2}, \dots, \ls_\nu$ are the $\nu$-admissible, ordered lists used in defining the boundary system $\B_\nu.$
\end{enumerate}
\end{stratcatlin}

\smallskip\noindent {\bf Remark:} The difference between this statement and Catlin's original Proposition 2.1 in \cite{catlinbdry} is in shrinking $U$ such that all properties hold on the closure of $U$ in $b \Omega,$ $\overline{U}.$ Parts (i)-(iii) hold for a given neighborhood as shown by Catlin, so they will clearly hold on any smaller neighborhood of $x_0.$ Catlin proved that properties in (iv) and (v) hold at $x_0.$ These are open conditions, however, and there are only finitely many of them, so it is obvious the neighborhood $U$ can be shrunk, if necessary, so that they hold on the closure $\overline{U}$ of the shrunken neighborhood. Furthermore, note that condition (v) implies that the gradients of the functions $r_1, r_{p+2}, \dots, r_\nu$ are nonzero on $\overline{U}$ and linearly independent, which makes $M^\nu$ a manifold as stated in (ii).

\smallskip The stratification induced by the commutator multitype is most interesting when there are only finitely many strata in a neighborhood of a point $x_0.$ Here are two important cases when the number of strata is finite: \label{finnumstrata}

\begin{enumerate}
\item[(i)] The domain $\Omega$ is pseudoconvex and of finite D'Angelo type at $x_0.$
\item[(ii)] The defining function $r$ is a polynomial.
\end{enumerate}

\noindent Case (ii) is easy to see as $\ct (x_0)$ is constructed by differentiation, and a polynomial only has finitely many non-zero derivatives. 

Case (i) follows from a number of the results stated above. Let us assume the D'Angelo $q$-type $\Delta_q (b \Omega, x_0)$ is finite. As D'Angelo proved in \cite{opendangelo}, which we stated as Theorem~\ref{openqtthm}, the finiteness of the D'Angelo type is an open condition and the D'Angelo $q$-type is locally bounded, so we can shrink the neighborhood $U$ from Proposition~\ref{stratcatlinprop} to a neighborhood $\tilde U,$ where $x_0 \in \tilde U$ and $\tilde U \subset U,$ so that the D'Angelo type is finite at all $x \in \tilde U\cap b \Omega$ and $\Delta_q (b \Omega, x) \leq 2 (\Delta_q(b \Omega, x_0))^{n-q}.$ Part (4) of Theorem~\ref{mlttheorem} guarantees that all entries up to the $(n+1-q)^{th}$ one of $\mt(x)$ are controlled by $\Delta_q (b \Omega, x)$ and $\mt(x) = \ct (x)$ for a pseudoconvex domain by part (ii) of Theorem~\ref{cmtypethm}. Furthermore, part (i) of Proposition~\ref{stratcatlinprop} guarantees the upper semi-continuity of the commutator multitype on $\tilde U\cap b \Omega.$ By the definition of the set of weights $\Gamma_n$ and the fact that $\mt (x),$ $\ct (x) \in \Gamma_n,$ the number of level sets of $\ct^{n+1-q}$ must be finite everywhere in $\tilde U \cap b \Omega.$ This observation appears as remark 1.2 on page 532 of Catlin's paper \cite{catlinbdry}. 

Let us now assume we are either in case (i) or case (ii). The $(n+1-q)^{th}$ commutator multitype $\ct^{n+1-q}$ takes only finitely many values $\ct^{n+1-q}_1, \dots, \ct^{n+1-q}_N$ at all points of $\tilde U \cap b \Omega,$ where $\ct^{n+1-q}_1< \ct^{n+1-q}_2 < \cdots < \ct^{n+1-q}_N$ and $N$ is some natural number, $N \geq 1.$ $\tilde U$ is given by Theorems~\ref{stratcatlinprop} and ~\ref{openqtthm} in case (i), whereas $\tilde U$ does not have to satisfy any conditions in case (ii). Let $$S_j= \{ x \in \tilde U \cap b \Omega \: \big| \: \ct^{n+1-q} (x)=\ct^{n+1-q}_j\}$$ be the level sets of the $(n+1-q)^{th}$ commutator multitype for $1 \leq j \leq N.$ Note that unlike in Proposition~\ref{stratcatlinprop}, we are working here with the open set $\tilde U$ and {\bf not} its closure. Clearly, $$\tilde U \cap b \Omega = \bigcup_{j=1}^N \, S_j \:\: {\text and} \:\: S_i \cap S_j = \emptyset \:\: {\text for} \:\: i \neq j$$ We shall now show that $S_1,$ the level set of the lowest commutator multitype, is an open set in the induced topology on $b \Omega.$

\medskip
\newtheorem{lowopen}[wgt]{Lemma}
\begin{lowopen}
\label{lowopenlemma}
Let $\Omega \subset \C^n$ be a smooth domain, and let $x_0 \in b \Omega$ have a neighborhood $\tilde U$ such that the $\nu^{th}$ commutator multitype assumes only finitely many values $\ct^\nu_1< \ct^\nu_2 < \cdots < \ct^\nu_N$ in $\tilde U \cap b \Omega$ for $N$ some natural number $N \geq 1$ and $\nu \geq 2.$ Let $S_j$ be the level set of points in $\tilde U \cap b \Omega$ with $\nu^{th}$ commutator multitype $\ct^\nu_j.$ The level set $S_1$  of the lowest commutator multitype is open in $b \Omega.$
\end{lowopen}

\smallskip\noindent {\bf Proof:}  $\ct^\nu_1= (1, c_2, \dots, c_\nu).$ We distinguish two cases:

\medskip\noindent {\bf Case 1:} $c_2 = \cdots = c_\nu = + \infty.$ Since $\ct^\nu_1$ is the lowest $\nu^{th}$ commutator type in $\tilde U \cap b \Omega,$ then every single point of $\tilde U \cap b \Omega$ has $\ct^\nu_1= (1, +\infty, \dots, +\infty).$ $S_1 = \tilde U \cap b \Omega$ and is thus open. \qed

\medskip\noindent {\bf Case 2:} There exists $c_\mu < +\infty$ among $c_2, \dots, c_\nu$ in $\ct^\nu_1.$ Let $\mu$ be the highest integer among $2, \dots, \nu$ for which this condition holds. For any $x \in S_1,$ the entries of $\ct^\mu(x)$ are finite, so there exists some $\B_\mu (x),$ which is a boundary system of codimension $n-\mu.$ By part (i) of Proposition~\ref{stratcatlinprop}, there exists a neighborhood $U_x \ni x$ such that for every $y \in U_x,$ $\ct^\mu (y) \leq \ct^\mu (x).$ It follows that $\ct^\mu (y) = \ct^\mu (x)$ because if $\ct^\mu (y) < \ct^\mu (x),$ then $\ct^\nu (y) < \ct^\nu (x),$ which is impossible since $x \in S_1,$ the lowest level set of the $\nu^{th}$ commutator multitype. If $\mu=\nu,$ we are done as $y \in S_1$ for every $y \in U_x.$ For $\mu <\nu,$ we are still able to conclude from $\ct^\mu (y) = \ct^\mu (x)$ that $\ct^\nu (y) = \ct^\nu (x)$ for all $y \in U_x$ since $\ct^\nu (y) < \ct^\nu (x)$ leads to the same contradiction as before whereas $\ct^\nu (y) > \ct^\nu (x)$ is impossible as $\ct^\nu (y)$ and $\ct^\nu (x)$ agree up to and including the $\mu^{th}$ entry and $c_{\mu+1}=\cdots = c_\nu = + \infty.$ Once again, $y \in S_1$ for every $y \in U_x.$ \qed

\bigskip\noindent When there is a neighborhood of a boundary point $x_0$ that contains only finitely many level sets of a partial commutator multitype $\ct^\nu$ {\bf and} all entries of $\ct^\nu$ are finite, then it turns out that $\ct^\nu_1,$ the lowest partial commutator multitype in that neighborhood, has the lowest possible value, $\ct^\nu_1 = (1, 2, \dots, 2).$  We shall prove this result for $\nu=n+1-q$ first in the case when $\Omega$ is a pseudoconvex smooth domain of finite D'Angelo $q$-type and then in the case when $\Omega$ is defined by a polynomial and finite D'Angelo $q$-type still holds. The latter will be crucial for the proof of the Main Theorem, Theorem~\ref{maintheorem}.

\medskip
\newtheorem{lowdef}[wgt]{Lemma}
\begin{lowdef}
\label{lowdeflemma} Let $\Omega \subset \C^n$ be a pseudoconvex smooth domain, and let $x_0 \in b \Omega$ be a boundary point of finite D'Angelo $q$-type. Let $\tilde U$ be a neighborhood of $x_0$ such that on $\tilde U \cap b \Omega,$ $b \Omega$ has finite D'Angelo $q$-type everywhere and the $(n+1-q)^{th}$ commutator multitype $\ct^{n+1-q}$ takes only finitely many values $\ct^{n+1-q}_1< \dots< \ct^{n+1-q}_N$ for some natural number $N \geq 1.$ The lowest $(n+1-q)^{th}$ commutator multitype $\ct^{n+1-q}_1=(1,2, \dots,2).$ 
\end{lowdef}

\smallskip\noindent {\bf Proof:}  $\ct^{n+1-q}_1= (1, c_2, \dots, c_{n+1-q}).$ Theorem~\ref{openqtthm}, part (4) of Theorem~\ref{mlttheorem}, and part (ii) of Theorem~\ref{cmtypethm} together imply that $c_2 \leq \dots \leq c_{n+1-q} < +\infty$ as explained above. Now assume there exists $c_\mu >2$ among $c_2, \dots, c_{n+1-q}$ in $\ct^{n+1-q}_1.$ Let $\mu$ be the highest integer among $2, \dots, n+1-q$ for which this condition holds. Consider $x \in S_1.$ All entries of $\ct^\mu$ are finite, so there exists some $\B_\mu (x),$ which is a boundary system of codimension $n-\mu$ at $x.$ Since $2 < c_\mu < +\infty,$ $\B_\mu$ contains a real valued function $r_\mu$ whose gradient is linearly independent from the gradient of the defining function $r.$ Parts (ii) and (iii) of Proposition~\ref{stratcatlinprop} show that there is a neighborhood $U_x \ni x$ such that $S_1 \cap U_x \subset \{y \in U_x \cap b \Omega   \: \big| \: r_\mu (y)=0\}.$ From the previous result, Lemma~\ref{lowopenlemma}, however, we know $S_1$ is open in $b \Omega,$ so it cannot be contained in a set of codimension $1$ in $b \Omega.$ Therefore, no such function $r_\mu$ can exist, and $\ct^{n+1-q}_1= (1, 2, \dots, 2).$\qed

\medskip
\newtheorem{lowpol}[wgt]{Lemma}
\begin{lowpol}
\label{lowpollemma} Let $\Omega \subset \C^n$ be a domain given by a defining function $r$ that is a polynomial, and let $x_0 \in b \Omega$ be a boundary point of finite D'Angelo $q$-type. Let $\tilde U$ be a neighborhood of $x_0$ such that on $\tilde U \cap b \Omega,$ $b \Omega$ has finite D'Angelo $q$-type everywhere and the $(n+1-q)^{th}$ commutator multitype $\ct^{n+1-q}$ takes only finitely many values $\ct^{n+1-q}_1< \dots< \ct^{n+1-q}_N$ for some natural number $N \geq 1.$ The lowest $(n+1-q)^{th}$ commutator multitype $\ct^{n+1-q}_1=(1,2, \dots,2).$ 
\end{lowpol}

\smallskip\noindent {\bf Proof:}  $\ct^{n+1-q}_1= (1, c_2, \dots, c_{n+1-q}).$ Assume there exists $c_\mu >2$ among $c_2, \dots, c_{n+1-q}$ in $\ct^{n+1-q}_1.$ We distinguish two cases:

\medskip\noindent {\bf Case 1:} $2 <c_\mu < +\infty.$ Let $x \in S_1.$ We argue as in the proof of Lemma~\ref{lowdeflemma} by looking at a boundary system $\B_\mu (x)$ at $x$ and at the real valued function $r_\mu$ corresponding to the entry $c_\mu.$ We arrive at a contradiction of the openness of the set $S_1.$ This case is thus impossible.

\medskip\noindent {\bf Case 2:} $c_\mu = +\infty.$ Let $\mu$ be the lowest integer among $2, \dots, n+1-q$ for which this condition holds. By the argument in Case 1, $\ct^{\mu-1}_1 = (1,2, \dots, 2)$ unless $\mu=2.$ Let $x \in S_1.$ By  Lemma~\ref{lowopenlemma}, we know there exists a neighborhood $U_x \ni x$ such that $U_x \cap b \Omega \subset S_1.$ All points $y$ of $U_x \cap b \Omega$ have $(n+1-q)^{th}$ commutator multitype $\ct^{n+1-q} (y) = (1, 2, \dots, 2, + \infty, \dots, +\infty)$ if $\mu >2$  or $\ct^{n+1-q} (y) = (1, + \infty, \dots, +\infty)$ otherwise. The full commutator multitype $\ct$ thus has at least q entries of $+\infty.$ An entry of $+\infty$ signals a direction that belongs to the null space of the Levi form.  Therefore, the dimension of the null space of the Levi form is constant at every $y \in U_x \cap b \Omega$ and satisfies $\dim \ny \geq q.$ The foliation result of Freeman and Sommer (\cite{sommer1}, \cite{sommer2}, \cite{freeman1}, \cite{freeman2}), which holds regardless of pseudoconvexity, implies the open set $U_x \cap b \Omega$ in the boundary of the domain is foliated by complex manifolds of dimension equal to this constant dimension of the null space of the Levi form, which is at least $q.$ This violates finite D'Angelo $q$-type that is assumed to hold on all of $\tilde U \cap b \Omega.$ We have shown this case is also impossible.

\medskip\noindent We conclude $c_2 = \dots = c_{n+1-q} =2.$ \qed

\smallskip\noindent {\bf Remark:} The proof of Lemma~\ref{lowpollemma} is complicated by the fact that in the absence of pseudoconvexity, entries of $+\infty$ cannot be ruled out from $\ct^{n+1-q}$ for two reasons: The equality between the commutator multitype $\ct$ and the multitype $\mt$ might no longer hold and also the domination of the entries of the multitype by the D'Angelo $q$-type could fail. The interested reader should take a look at Theorem 3.7 on p.543 of \cite{catlinbdry} where David Catlin shows $\ct (x_0) \geq \mt(x_0)$ even in the absence of psedoconvexity. Unfortunately, the inequality in the other direction is proven using pseudoconvexity in a rather fundamental way. As for the other result, as can be seen on p.556 of \cite{catlinbdry}, the domination of entries of $\mt(x_0)$ by the D'Angelo $q$-type $\Delta_q (b \Omega, x_0)$ requires the existence of a particular system of coordinates, whose existence Catlin proves using normalization, a technique that requires pseudoconvexity. Perhaps Martin Kol\'a\v{r}'s recent work quoted in the introduction could shed some light on this point in the absence of pseudoconvexity. The upper semi-continuity of $\mt$ is also plausible but not evident when the domain is not pseudoconvex. The main obstacle is that $\mt$ has a very abstract definition. When the domain is pseudoconvex, Catlin can characterize $\mt$ by relating it to $\ct,$ which can be explicitly computed. When the equality of $\mt$ and $\ct$ is not known to hold, clearly more work is necessary.

\medskip\noindent We end this section with a corollary to the previous result:

\medskip
\newtheorem{polnonvan}[wgt]{Corollary}
\begin{polnonvan}
\label{polnonvancor} Let $\Omega \subset \C^n$ be a domain given by a defining function $r$ that is a polynomial, and let $x_0 \in b \Omega$ be a boundary point of finite D'Angelo $q$-type. The Levi determinant $\text{coeff}\{\partial r \wedge \bar \partial r \wedge (\partial \bar \partial r)^{n-q}\}$ does not vanish to infinite order.
\end{polnonvan}

\smallskip\noindent {\bf Proof:} This corollary concerns case (ii) on page~\pageref{finnumstrata} from the discussion of instances when the commutator multitype defines finitely many strata. By Theorem~\ref{openqtthm}, there exists a neighborhood $\tilde U \ni x_0$ such that on $\tilde U \cap b \Omega,$ $b \Omega$ has finite D'Angelo $q$-type everywhere. All hypotheses of the previous result, Lemma~\ref{lowpollemma}, are now satisfied, so we can now apply it to conclude that the lowest achieved value of the $(n+1-q)^{th}$ commutator multitype $\ct^{n+1-q}$ has to be $\ct^{n+1-q}_1=(1,2, \dots,2)$ in $\tilde U \cap b \Omega.$ Each entry of $2$ in the commutator multitype indicates a non-zero Levi eigenvalue. In other words, at every point where $\ct^{n+1-q}_1=(1,2, \dots,2)$ is achieved, there are $n-q$ non-zero Levi eigenvalues, so the Levi determinant $\text{coeff}\{\partial r \wedge \bar \partial r \wedge (\partial \bar \partial r)^{n-q}\}$ is non-zero. A polynomial that is non-zero at at least one point cannot vanish to infinite order, so $\text{coeff}\{\partial r \wedge \bar \partial r \wedge (\partial \bar \partial r)^{n-q}\}$ does not vanish to infinite order. \qed

\section{Proof of the Main Theorem}
\label{mainthmpf}

We will start this section with a crucial example that sheds light on the hypotheses of the main theorem, Theorem~\ref{maintheorem}. At the end of last section, we proved Corollary~\ref{polnonvancor} whose conclusion was that for a domain defined by a polynomial $r$ whose boundary has at least one point of finite D'Angelo $q$-type, the Levi determinant $\text{coeff}\{\partial r \wedge \bar \partial r \wedge (\partial \bar \partial r)^{n-q}\}$ cannot vanish to infinite order. Since Corollary~\ref{polnonvancor} is based on Lemma~\ref{lowpollemma} whose proof involves understanding the stratification given by levels of the Catlin commutator multitype, the reader might be lead to believe the condition of finite D'Angelo $q$-type can be replaced by simply requiring finite entries of $\ct^{n+1-q}$ at $x_0.$ This is unfortunately false as the following example drawn from p.217 of \cite{kohndangelo} shows: $$r=2 \, Re \{ z_1 \} + |z_2^2-z^3_3|^2.$$ This is a domain in $\C^3$ defined by a polynomial. Let $q=1.$ $\ct^2 (0)= (1, 4,6),$ so all entries of the commutator multitype are finite at the origin, but D'Angelo $1$-type fails to be finite there as the holomorphic curve $\varphi(t)=(0, t^3, t^2)$ for $t \in \C$ sits in $b \Omega$ and passes through the origin. In fact, $\ct^2 (x)= (1, 2,+\infty)$ away from the origin as the Levi form has rank $0$ at the origin and $1$ everywhere else. The determinant of the Levi form $\text{coeff}\{\partial r \wedge \bar \partial r \wedge (\partial \bar \partial r)^2\}$ can easily be shown to be identically zero in this case, so the conclusion of Lemma~\ref{lowpollemma} fails.

We can now proceed with the proof of the main result of this paper:

\medskip\noindent {\bf Proof of Theorem~\ref{maintheorem}:} $\Omega \subset \C^n$ is a domain defined by a smooth real-valued function $r.$ We are assuming the D'Angelo $q$-type at $x_0 \in b \Omega$ is finite, i.e. $\Delta_q (b \Omega, x_0)=t.$ Consider the Taylor expansion of $r$ at $x_0.$ Since the D'Angelo $q$-type is a finitely determined condition, the remark after Theorem~\ref{findetqtthm} implies that any truncation $\tilde r$ of the Taylor expansion of $r$ at $x_0$ such that the following two conditions hold:

\begin{enumerate}
\item $\tilde r$ is real valued;
\item $\tilde r$ contains all terms of the Taylor expansion of $r$ at $x_0$ that are of order $\lceil t \rceil,$ the roundup of $t,$ and lower;
\end{enumerate}
satisfies that if $\tilde \Omega$ is the domain defined by $\tilde r,$ then the D'Angelo $q$-type of $b \tilde \Omega$ at $x_0$ is the same, namely $\Delta_q (b \tilde \Omega, x_0)=t.$ We apply Corollary~\ref{polnonvancor} to the domain defined by polynomial $\tilde r$ to conclude $\text{coeff}\{\partial \tilde r \wedge \dbar \tilde r \wedge (\partial \dbar \tilde r)^{n-q}\}$ cannot vanish to infinite order. If $\tilde r$ is a polynomial of degree exactly $\lceil t \rceil,$ then the order of vanishing of its Levi determinant $\text{coeff}\{\partial \tilde r \wedge \dbar \tilde r \wedge (\partial \dbar \tilde r)^{n-q}\}$ can be at most $(\lceil t \rceil -2 )^{n-q}.$ Let us show the same is true of the Levi determinant of the original defining function $r.$

We first prove that the order of vanishing of the Levi determinant $\text{coeff}\{\partial \tilde r \wedge \dbar \tilde r \wedge (\partial \dbar \tilde r)^{n-q}\}$ of any truncation $\tilde r$ of the Taylor expansion of $r$ at $x_0$ satisfying conditions (1) and (2) is at most $(\lceil t \rceil -2 )^{n-q}.$ Assume not, i.e. assume that there exists a truncation $\tilde r$ such that the order of vanishing of 
$\text{coeff}\{\partial \tilde r \wedge \dbar \tilde r \wedge (\partial \dbar \tilde r)^{n-q}\}$ is at least $(\lceil t \rceil -2 )^{n-q} +1.$ Now truncate $\tilde r$ by throwing out all terms of degree strictly greater than $\lceil t \rceil.$ We have a polynomial $\tilde{ \tilde r}$ that satisfies conditions (1) and (2) along with $\Delta_q (b \tilde {\tilde \Omega}, x_0)=t,$ where $\tilde {\tilde \Omega}$ is the domain defined by $\tilde{ \tilde r},$ and $\text{coeff}\{\partial \tilde{ \tilde r}\wedge \dbar \tilde{ \tilde r} \wedge (\partial \dbar \tilde{ \tilde r})^{n-q}\}$ vanishes to infinite order by construction contradicting Corollary~\ref{polnonvancor}.

Finally, since every truncation of the Taylor expansion of $r$ at $x_0$ satisfying conditions (1) and (2) has a Levi determinant $\text{coeff}\{\partial \tilde r \wedge \dbar \tilde r \wedge (\partial \dbar \tilde r)^{n-q}\}$ that vanishes to order at most $(\lceil t \rceil -2 )^{n-q},$ then the same must be true of the original defining function $r.$ \qed

\bibliographystyle{plain}
\bibliography{SmoothTypeEquiv}

\end{document}